\documentclass[12pt]{article}
\usepackage{amssymb,amsmath,amscd}
\usepackage{dsfont}
\usepackage{verbatim}
\usepackage{amsmath} 
\usepackage{pstricks}
\usepackage{enumitem}
\usepackage{ntheorem}
\usepackage{amssymb}
\usepackage{breqn}
\usepackage[left=2cm,top=2cm,right=2cm,bottom=2cm,nohead]{geometry}

\newtheorem{de}{Definition}[section]

\newtheorem{theo}{Theorem}

\newtheorem{conj}{Conjecture}[section]

\setenumerate[1]{font=\bfseries,label=\Roman*.}
\setenumerate[2]{font=\itshape,label=\arabic*)}

\author{Yashar Memarian
}       

\title{Spherical Localisation in Convex and Metric Geometry \\
       }
\date{\today}

\begin{document}
\maketitle
\begin{abstract}
Spherical localisation is a technique whose history goes back to M.Gromov and V.Milman. It's counterpart, the Euclidean localisation is extensively studied and has been put to great use in various branches of mathematics. The purpose of this paper is to quickly introduce spherical localisation, as well as demonstrate some of its applications in convex and metric geometry.
\end{abstract}
\subsection{Spherical Localisation}

Spherical localisation is the spherical counterpart of what is known as Euclidean localisation. It was first used in \cite{gromil}. The canonical sphere $\mathbb{S}^n$ is the Riemannian sphere with sectional curvature equal to $1$ and its Riemannian metric structure. For simplicity, we view this space as a probability space upon which we have normalised the Riemannian volume. We shall denote the normalised measure by $d\mu$. First we need to understand what the counterpart of a \emph{needle} on the sphere is.
\begin{de}[Spherical Needle]
A  $k$-dimensional spherical needle is a pair $(S,\nu)$, where $S$ is a $k$-dimensional convex subset of $\mathbb{S}^n$ and $\nu$ is a probability measure defined on $S$ such that $d\nu(x)=h(x) dv(x)$. Where $h$ is a function (defined on $S$) such that its restriction on every geodesic segment inside $S$ can be written as $C\sin(t+t_0)^{n-k}$, where $C>0$ and $t_0\in[0,\pi]$ and $t$ is the parametrisation of the geodesic segment by its arc length, and $dv$ is the $k$-dimensional Riemannian volume assigned to the $k$-dimensional sphere $\mathbb{S}^k$ which contains the convex set $S$. When $k=1$, we refer to $S$ as a $1$-needle.
\end{de}

The task of the spherical localisation can be viewed in the following theorem which, is proved in \cite{memwst} and \cite{memgauss}:

\begin{theo} \label{local}
Let $k<n$. Let $f_1, f_2,\cdots, f_{k+1}$ be $k+1$ continuous functions on $\mathbb{S}^n$. Suppose that for every $i\in\{1,\cdots,k+1\}$ we have:
\begin{eqnarray*}
\displaystyle\int_{\mathbb{S}^n}f_i(u)d\mu(u)>0.
\end{eqnarray*}
Then, there exists a (measurable convex) partition of $\mathbb{S}^n$ into at most $k$-dimensional needles such that for every $(S,\nu)$ needle in this partition and every $i\in\{1,\cdots,k+1\}$ we have:
\begin{eqnarray*}
\displaystyle\int_{S}f_i(u)d\nu(u)>0.
\end{eqnarray*}
\end{theo}

For a rigorous definition of a measurable convex partition of the sphere we invite the reader to consult \cite{memwst}. Even if Theorem \ref{local} and spherical localisation seem like \emph{straightforward} generalisations of Euclidean localisation, we see in the following sections that they impart very important results in different branches of geometry.

\section{Spherical Localisation in Metric Geometry}

The \emph{first} serious application of spherical localisation is in Gromov's famous waist of the sphere Theorem :
\begin{theo}[Waist of the Sphere Theorem] \label{a}
Let $f:\mathbb{S}^n\to\mathbb{R}^k$ be a continous map. There exists a point $z\in\mathbb{R}^k$ such that for every $\varepsilon>0$ we have:
\begin{eqnarray*}
vol_n(f^{-1}(z)+\varepsilon)\geq vol_n(\mathbb{S}^{n-k}+\varepsilon),
\end{eqnarray*}
where $vol_n(X+\varepsilon)$ denotes the volume of the $\varepsilon$-neighborhood of the set $X$.
\end{theo}

For the proof, one can consult the original paper of Gromov \cite{grwst} or the more recent paper of the author \cite{memwst}. Later on, in \cite{roman}, the authors generalise the above result in which one consider maps between spheres to general $k$-dimensional manifolds. In \cite{memusphere}, the author, once again applying spherical localisation, proves a waist theorem for unit spheres of uniformly-convex normed spaces :

\begin{theo}[Waist of Unit Spheres of Uniformly-Convex Normed Spaces] \label{b}
Let $X$ be a uniformly convex normed space of finite dimension $n+1$. Let $S(X)$ be the unit sphere of $X$, for which the distance is induced from the norm of $X$. The measure defined on $S(X)$ is the conical probability measure denoted by $\mu$. Then for every continuous map $f:S(X)\to\mathbb{R}^{k}$ there exists a $z\in\mathbb{R}^k$ such that for every $\varepsilon>0$ we have:
\begin{eqnarray*}
\mu(f^{-1}(z)+\varepsilon)\geq w(\varepsilon),
\end{eqnarray*}
where
\begin{eqnarray*}
w(\varepsilon)= \frac{1}{1+(1-2\delta(\frac{\varepsilon}{2}))^{n-k}(k+1)^{k+1}\frac{F(k,\frac{\varepsilon}{2})}{G(k,\frac{\varepsilon}{2})}}
\end{eqnarray*}
where $\delta(\varepsilon)$ is the modulus of convexity ,
\begin{eqnarray*}
F(k,\varepsilon)= \int_{\psi_2(\varepsilon)}^{\frac{\pi}{2}}\sin(x)^{k-1}\,dx.
\end{eqnarray*} 
and 
\begin{eqnarray*}
G(k,\varepsilon)= \int_{0}^{\psi_1(\varepsilon)}\sin(x)^{k-1}\,dx.
\end{eqnarray*}
And where 
\begin{eqnarray*}
\psi_1(\varepsilon)=2\arcsin(\frac{\varepsilon}{4\sqrt{k+1}})
\end{eqnarray*}
and
\begin{eqnarray*}
\psi_2(\varepsilon)=2\arcsin(\frac{\varepsilon}{2\sqrt{k+1}})
\end{eqnarray*}
\end{theo}

Both theorems \ref{a} and \ref{b} are proved using suitable partitions of the \emph{spheres} into needles and translating the global problem into smaller needles where the computations are much easier.

\section{Spherical Localisation in Convex Geometry}

Spherical Localisation has very serious applications in convex geometry. So far, by applying this technique, one is finally able to answer the old-standing Gaussian Correlation Conjecture, as well as also provide an interesting lower bound for the Mahler volume of symmetric convex bodies. We review these two important results here:
\subsection{The Gaussian Correlation Conjecture}

The Gaussian Correlation Conjecture is the following:
\begin{conj}\label{1}
Let $\gamma_n$ be the Gaussian probability measure on $\mathbb{R}^n$. Given two symmetric convex bodies $K_1$ and $K_2$ in $\mathbb{R}^n$, the following holds:
\begin{equation} \label{eqn:ti}
\gamma_n(\mathbb{R}^n)\gamma_n(K_2\cap K_2)\geq \gamma_n(K_1)\gamma_n(K_2).
\end{equation}
\end{conj}

Since this conjecture has occupied mathematicians for many years, it is useful to provide a sketch of the proof here in order to demonstrate how spherical localisation can be used to prove it.
Let us work a bit more on the inequality (\ref{eqn:ti}). Clearly, inequality (\ref{eqn:ti}) is an inequality involving four continuous functions defined on the canonical sphere $\mathbb{S}^{n-1}$. Indeed, let us express the Gaussian measure defined on $\mathbb{R}^n$ in polar coordinates associated to $\mathbb{R}^n$. In this case, we shall have:
\begin{eqnarray*} 
\gamma_n=r^{n-1}e^{-r^2/2}dr\,du_{\mathbb{S}^{n-1}},
\end{eqnarray*}
where $du_{\mathbb{S}^{n-1}}$ is the canonical Riemannian measure associated to the canonical sphere. Therefore inequality \ref{eqn:ti} is written in the following form:
\begin{equation} \label{eqn:ti2}
\int_{\mathbb{S}^{n-1}}f_1(u)du \int_{\mathbb{S}^{n-1}}f_2(u)du\geq \int_{\mathbb{S}^{n-1}}f_3(u)du \int_{\mathbb{S}^{n-1}}f_4(u) du,
\end{equation}
where for $i=1,2,3,4$, $f_i$ is a continuous function defined on $\mathbb{S}^{n-1}$ defined as follows. Set $M_1=\mathbb{R}^n$, $M_2=K_1\cap K_2$, $M_3=K_1$ and $M_4=K_2$. Therefore the function $f_i$ takes a point $u\in\mathbb{S}^{n-1}$ and calculates the following :
\begin{eqnarray*}
f_i(u)=\int_{0}^{x_i(u)}r^{n-1}e^{-r^2/2}dr,
\end{eqnarray*}
where $x_i(u)$ is the the length of the segment starting from $0$, going in the direction $u$, and touching the boundary of $M_i$. Of course for $M_1$, for every $u\in\mathbb{S}^{n-1}$, we have $x_1(u)=+\infty$. 

Applying (a mild version of) Theorem \ref{local}, one can illustrate that in order for inequality \ref{eqn:ti2} to be true, it is enough that for every pair $(\sigma,\nu)$ of $1$-needle, the following inequality will be true :
\begin{equation} \label{eqn:ti3}
\int_{\sigma}f_1(t)d\nu(t)\int_{\sigma}f_2(t)d\nu(t)\geq \int_{\sigma}f_3(t)d\nu(t)\int_{\sigma}f_4(t)d\nu(t).
\end{equation}
Inequality (\ref{eqn:ti3}) tells us to consider the restriction of these functions on $\sigma$. Taking the restriction of the functions $f_i$ on $\sigma$ means that we should intersect the sets $\mathbb{R}^n$, $K_1$, $K_2$ and $K_1\cap K_2$ with a $2$-dimensional plane which contains the origin of $\mathbb{R}^n$ and the geodesic segment $\sigma$. On this $2$-dimensional plane, we consider the cone which is defined over $\sigma$. It means this cone contains the origin of $\mathbb{R}^n$ and one half-line on the boundary of this cone starts from the origin, touches one point of the boundary of $\sigma$, and goes to infinity. The same for the other half-line. We denote this cone by $C$. Now we can write down the inequality (\ref{eqn:ti3}) on $C$. Parametrising $C$ in polar coordinates, we will obtain the measure defined on $C$ which will be given by $\mu_2=r^{n-1}e^{-r^2/2}\cos(t+t_0)^{n-2}dr\,dt$. Therefore, inequality (\ref{eqn:ti3}) is equivalent to the following:
\begin{equation} \label{eqn:ti4}
\mu_2(C\cap \mathbb{R}^2)\mu_2(K_1\cap K_2\cap C)\geq \mu_2(K_1\cap C)\mu_2(K_2\cap C).
\end{equation}

We are then able to simplify the verification of Conjecture \ref{1} to the verifications of a family of $2$-dimensional problems given by inequality (\ref{eqn:ti4}). This means that if we can prove that for every cone $C$ and every measure $\mu_2$ of the form given above, if inequality (\ref{eqn:ti4}) is true then Conjecture \ref{1} is true.


There is a nice symmetrisation method which works well (only in dimension $2$) for studying problems such as the one given by inequality (\ref{eqn:ti4}). We can show that for every $C$, $K_1$, $K_2$ and $\mu_2$ there are always two symmetric strips $S_1$ and $S_2$ such that :
\begin{eqnarray*}
\frac{\mu_2(C\cap \mathbb{R}^2)\mu_2(K_1\cap K_2\cap C)}{\mu_2(K_1\cap C)\mu_2(K_2\cap C)}\geq \frac{\mu_2(\mathbb{R}^2)\mu_2(S_1\cap S_2)}{\mu_2(S_1)\mu_2(S_2)}.
\end{eqnarray*}

Now if for every pair of symmetric strips $S_1$ and $S_2$  and every measure $\mu_2$ we have:
\begin{equation} \label{eqn:ti5}
\mu_2(\mathbb{R}^2)\mu_2(S_1\cap S_2)\geq \mu_2(S_1)\mu_2(S_2),
\end{equation}
then Conjecture \ref{1} would be true!

But unfortunately this is not true.

To prove Conjecture \ref{1}, we need to work harder. 

The trick is that the \emph{partition} version of spherical localisation (Theorem \ref{local}) allows one to \emph{not} have to consider all the $1$-needles. Let us explain this here. In the definition of the measure $\nu$, we had to use a parameter $t_0$. Remember $\nu=\cos(t+t_0)^{n-2}dt$. It is the polar part of the measure $\mu_2$. What is always true is the following : For every $C$ there is always a $t_0$ such that the following is true :
\begin{equation} \label{eqn:ti6}
\mu_{2,t_0}(C\cap \mathbb{R}^2)\mu_{2,t_0}(K_1\cap K_2\cap C)\geq \mu_{2,t_0}(K_1\cap C)\mu_{2,t_0}(K_2\cap C).
\end{equation}
We include the dependence on $t_0$ for the measure $\mu_2$ by writing it as $\mu_{2,t_0}$. 


Now assume Conjecture \ref{1} is not true. This means inequality (\ref{eqn:ti2}) is not true. Applying Theorem \ref{local}, if (\ref{eqn:ti2}) is not true, then there exists a \emph{partition} of $\mathbb{S}^{n-1}$ into geodesic needles, such that for every needle $(\sigma,\nu)$ in this partition the inequalities (\ref{eqn:ti3}), and hence (\ref{eqn:ti4}), are not true. By cleverly \emph{constructing} the partition of $\mathbb{S}^{n-1}$ into geodesic needles upon which inequalities (\ref{eqn:ti3}) and (\ref{eqn:ti4}) are not true, we can always find at least one geodesic segment $\sigma$ and a measure $\nu$ defined on this $\sigma$, such that translating the inequality on the cone $C$ defined over this $\sigma$, provides us with a measure $\mu_{2,t_0}$ (for which we know that inequality (\ref{eqn:ti6}) is \emph{always} true). This gives us a contradiction and finalises the demonstration of Conjecture \ref{1}.   

In the first draft and attempt of the proof of Conjecture \ref{1}, the author made two mistakes:
\begin{itemize}
\item The author only considered $1$-needles which have length equal to $\pi$. Indeed one has to consider $1$-needles which can very well have length strictly smaller than $\pi$.
\item The author attempted to prove inequality (\ref{eqn:ti3}) for \emph{every} $1$-needle. This, as it was said earlier, is false. To correct this mistake one needs the partition version of the spherical localisation, and has to show that by cleverly constructing the partition, and, by contradiction assuming the conjecture to be false, we obtain at least a spherical needle in the partition for which inequality (\ref{eqn:ti6}) is verified.
\end{itemize}

\subsection{Mahler Conjecture}
Similar to the Gaussian Correlation Conjecture, spherical localisation can also be used for studying the Mahler Conjecture. The Mahler conjecture is the following:
\begin{conj}
For every symmetric convex set in $\mathbb{R}^n$ where $n\geq 2$, we have:
\begin{eqnarray*}
vol(K)vol(K^{\circ})\geq \frac{4^n}{n!},
\end{eqnarray*}
where $K^{\circ}$ is the polar (dual) of $K$.
\end{conj}

In \cite{memmahler}, the author, in a manner similar to that used to prove Conjecture \ref{1}, proves the following:

\begin{theo} \label{main}
Let $K$ be a symmetric convex set in $\mathbb{R}^n$ where $n\geq 4$. Then
\begin{eqnarray*}
vol(K)vol(K^{\circ})&\geq& \alpha(n-1)vol_{n-1}(\mathbb{S}^{n-1})^2\\
                    &=&\frac{4\alpha(n-1)\pi^n}{\Gamma(\frac{n}{2})^2},
\end{eqnarray*}
where $\alpha(n)$ is defined as follows. Let $I\subset \mathbb{S}^1$ be a connected interval of the circle. Let $C(I)$ be the cone in $\mathbb{R}^2$ over $I$. Let $S(I)$ be the class of all connected closed intervals of the unit circle. Let $S$ be a $2$-dimensional convex set symmetric with respect to the origin of $\mathbb{R}^2$ containing the unit ball in its interior and being contained in the ball of radius $\sqrt{n+1}$. Define
\begin{eqnarray*}
\alpha(n,\theta,I,S)=\frac{\mu_{2,\theta}(C(I)\cap S)\mu_{2,\theta}(C(I)\cap S^{\circ})}{\big(\int_{I}g(\theta,t)dt\big)^2},
\end{eqnarray*}
where $\mu_{2,\theta}=r^n \cos(t+\theta)^{n-1}dr\wedge dt$  and $g(\theta,t)=\cos(t+\theta)^{n-1}$ is the density of this measure in polar coordinates. Define
\begin{eqnarray*}
\alpha(n,S)=\max_{\theta\in[0,\pi]}\big(\min_{I\in S(I)}\alpha(n,\theta,I,S)\big).
\end{eqnarray*}
And at last, define:
\begin{eqnarray*}
\alpha(n)=\min_{S}\big(\alpha(n,S)\big),
\end{eqnarray*}
where $S$ runs over the class of $2$-dimensional symmetric convex sets containing the unit ball in their interiors and being contained in the ball of radius $\sqrt{n+1}$. 
\end{theo} 

Since spherical localisation allowed one to give a \emph{sharp} result for the Gaussian Correlation conjecture, it is very interesting to see if the bound in Theorem \ref{main} is also \emph{sharp} or not.

All the problems in classical Convex Geometry which involve integral inequalities can be studied with spherical localisation technique. We examined two very important examples, but one can perhaps find many other problems of this kind.

\section{Spherical Localisation- an Overture to Algebraic Localisation?}
The sphere is one of the \emph{nicest} algebraic varieties. Is there a general \emph{theory} of localisation which can work for more general Algebraic varieties (maybe coming with Algebraic stratifications of Algebraic varieties)? Is what we have seen in the previous sections merely the tip of the iceberg for a vast machinery that could be utilised to study countless (Convex and) algebraic geometry problems? 

It's worth attempting!

\bibliographystyle{plain}
\bibliography{spherelo}

\begin{thebibliography}{1}

\bibitem{grwst}
M.~{Gromov}.
\newblock Isoperimetry of waists and concentration of maps.
\newblock {\em GAFA}, 13:178--215, 2003.

\bibitem{gromil}
M.~{Gromov} and V.D. {Milman}.
\newblock Generalisation of the spherical isoperimetric inequality to uniformly
  convex {B}anach spaces.
\newblock {\em Compositio Math.}, 62:3:263--282, 1987.

\bibitem{roman}
R.N.. {Karasev} and A.Yu. {Volovikov}.
\newblock Waist of the sphere for maps to manifolds.
\newblock {\em Topology and its Applications}, 160:13, 1987.

\bibitem{memwst}
Y.~{Memarian}.
\newblock On {G}romov's waist of the sphere theorem.
\newblock {\em Journal of {T}opology and {A}nalysis}, 3:7--36, 2011.

\bibitem{memusphere}
Y.~{Memarian}.
\newblock A lower bound on the waist of unit spheres of uniformly normed
  spaces.
\newblock {\em Compositio Math.}, 148(4):1238--1264, 2012.

\bibitem{memgauss}
Y.~Memarian.
\newblock The {G}aussian {C}orrelation {C}onjecture {P}roof.
\newblock {\em arXiv:1310.8099}, 2014.

\bibitem{memmahler}
Y.~Memarian.
\newblock A lower bound for the mahler volume of at least four-dimensional
  symmetric convex sets.
\newblock {\em arXiv:1501.02009}, 2015.

\end{thebibliography}

\end{document}